\documentclass[12pt]{article}
\usepackage{amssymb,latexsym}
\usepackage{amsbsy} 
\newtheorem{theo}{Theorem}

\newtheorem{lemm}[theo]{Lemma}

\newtheorem{prop}[theo]{Proposition}

\def\lb{[\![}
\def\rb{]\!]}

\def\ch{\mathop{\rm ch}\nolimits}
\def\mult{\mathop{\rm mult}\nolimits}
\def\rank{\mathop{\rm rank}\nolimits}
\def\tr{\mathop{\rm tr}\nolimits}
\def\al{\alpha}
\def\be{\beta}
\def\ga{\gamma}
\def\de{\delta}
\def\ep{\epsilon}

\def\th{\theta}
\def\la{\lambda}
\def\rh{\rho}
\def\si{\sigma}

\def\De{\Delta}

\def\C{\mathbb{C}}
\def\Z{\mathbb{Z}}
\def\beq{\begin{equation}}
\def\eeq{\end{equation}}
\def\bea{\begin{eqnarray}}
\def\eea{\end{eqnarray}}
\def\beas{\begin{eqnarray*}}
\def\eeas{\end{eqnarray*}}
\def\nn{\nonumber}

\def\mybox{\hfill\llap{$\Box$}}
\begin{document}
\addtolength{\baselineskip}{1mm}
\begin{center}
{\large \bf Fock representations of the Lie superalgebra $q(n+1)$}\\[1cm]
T Palev\footnote{Permanent address~: Institute for Nuclear Research
and Nuclear Energy, 1784 Sofia, Bulgaria 
(e-mail~: tpalev@inrne.bas.bg).} 
and J Van der Jeugt\footnote{Research Associate 
of the Fund for Scientific Research -- Flanders (Belgium);\\
(e-mail~: Joris.VanderJeugt@rug.ac.be).} \\
Department of Applied Mathematics and Computer Science,\\
University of Ghent, Krijgslaan 281-S9, B-9000 Gent, Belgium
\end{center}

\vskip 1cm

\begin{abstract}
For the Lie superalgebra $q(n+1)$ a description is given
in terms of creation and annihilation operators, in such a way
that the defining relations of $q(n+1)$ are determined by quadratic
and triple supercommutation relations of these operators.
Fock space representations $V_p$ of $q(n+1)$ are defined by means
of these creation and annihilation operators. These new
representations are introduced as quotient modules of some
induced module of $q(n+1)$. The representations $V_p$ are
not graded, but they possess a number of properties that
are of importance for physical applications.
For $p$ a positive integer, these
representations $V_p$ are finite-dimensional, with a unique highest
weight (of multiplicity 1). The Hermitian form that is consistent
with the natural adjoint operation on $q(n+1)$ is shown to be
positive definite on $V_p$.
For $q(2)$ these representations
are ``dispin''. For the general case of $q(n+1)$, 
many structural properties of $V_p$ are derived.
\end{abstract}

\newpage
%
\section{Introduction}

Lie superalgebras and their irreducible representations (simple
modules) have been the subject of much attention in both the
mathematical~\cite{Kac1,Kac2,Scheunert}
and the physics~\cite{Corwin,Balantekin,Hurni} literature.
However, even for the simplest family of basic classical Lie
superalgebras, namely $sl(m/n)$ or $gl(m/n)$, the understanding
of all finite-dimensional simple modules has been a
very difficult problem. The main reason for this difficulty
has been the existence of so-called {\em atypical} modules~\cite{Kac2}.
Although partial progress was made in determining a character
formula for atypical modules~\cite{VHKT1,VHKT2,KW}, 
the problem of determining the character for $gl(m/n)$ modules
was solved only recently~\cite{Serganova1,Serganova2} (see
also~\cite{VZ} for a simpler algorithm). 

The general linear Lie algebra $gl(n)$ has two super-analogues.
The first is the Lie superalgebra $gl(n/1)$, for which
the representations (even the atypical ones) are now well known~:
e.g., a Gelfand-Zetlin basis has been introduced and its
transformations have been determined~\cite{pal78,pal82},
and representations of $sl(n/1)$ have been studied~\cite{pal79,pal81}.
The second super-analogue is the strange Lie superalgebra $q(n)$.
Also this Lie superalgebra has received attention recently.
In particular, characters of finite-dimensional irreducible
graded representations of $q(n)$ have been 
determined~\cite{penkov0,penkov1,penkov2}, both in the
typical and atypical case. In a different context,
oscillator realizations have been given~\cite{frappat}.

The finite-dimensional irreducible graded representations of
$q(n)$ have the strange property that the multiplicity of the
highest (or lowest) weight is in general greater than 
1~\cite{penkov1}. From the physical point of view, where one
wishes to interpret the representation space as a Hilbert space
with a unique vacuum, this situation is not very favourable.

In the present paper, the purpose is to study a new class of
irreducible finite-dimensional
representations of the Lie superalgebra $q(n)$, which have 
certain properties that are required in a physical context.
In particular, the highest weight has multiplicity
1 (so there is a unique highest weight vector, up to a factor),
and the representation space can be naturally
equipped with a symmetric and 
nondegenerate positive definite Hermitian form (inner product).
Moreover, creation and annihilation 
generators (or operators) are introduced for
$q(n)$ such that the representation space is a Fock space.
The only property that has to be abandoned is the grading
of the representation space (but from the physical point of view,
this grading is no requirement).

The structure of the paper is as follows. In section~2 the main
definitions are given concerning the Lie superalgebra $q(n+1)$.
In section~3 a new basis for $q(n+1)$ is given in terms of
so-called creation and annihilation operators. The new class
of representations of $q(n+1)$ is introduced in section~4. These
representations are defined by means of an induced module
$\bar V_p$; the actual irreducible representation $V_p$ which
is of interest is then a quotient module of $\bar V_p$.
To have some idea about the structure of $V_p$, we first consider
the low rank case of $q(2)$ in section~5. Here, the representation
has a ``dispin'' structure. Section~6 goes back to the general
case $q(n+1)$, and includes several (technical) properties 
concerning the structure of $\bar V_p$, paving the way to
determining the structure of the simple modules $V_p$.
This is performed in section~7, where in particular we give
the dimension and character of $V_p$, and show that the Hermitian
form is positive definite.

\section{The Lie superalgebra $q(n+1)$}  

The Lie superalgebra $q(n+1)$ can be determined through its defining
representation, i.e.\
\beq
q(n+1) = \left\{ 
\left( \begin{array}{cc} A&B\\B&A \end{array} \right) \; | \;
A,B\in gl(n+1) \right\},
\eeq
where the matrices with $B=0$ are even, 
or elements of $q(n+1)_{\bar 0}$, and those
with $A=0$ odd, or elements of $q(n+1)_{\bar 1}$.
The subalgebra $sq(n+1)$ consists of those elements with $\tr(B)=0$. 
The Lie superalgebras $q(n+1)$ and $sq(n+1)$ still contain a 
one-dimensional center $\C I$, where $I$ is the identity matrix.
Hence one defines the quotient Lie superalgebra $psq(n+1)$ as 
$sq(n+1)/ \C I$. The notation $q(n+1)$, $sq(n+1)$ and $psq(n+1)$ is 
due to Penkov and others~\cite{penkov0}; 
in the notation of Kac~\cite{Kac1,Kac2} we have
$\bar Q (n)=sq(n+1)$ and $Q(n)=psq(n+1)$. 
Recall that $Q(n)$ is a simple Lie superalgebra for $n\geq 2$, and that
it is one of the series of classical (but ``strange'')
Lie superalgebras
in the classification of Kac~\cite{Kac1}.
For the development of representation theory, we shall be working mainly
with $q(n+1)$.

Thus $q(n+1)$ can be defined as the Lie superalgebra with $(n+1)^2$
even basis elements $e_{ij}^{\bar 0}$ ($i,j=0,1,\ldots,n$) and $(n+1)^2$
odd basis elements $e_{ij}^{\bar 1}$ ($i,j=0,1,\ldots,n$), 
satisfying the bracket relation
\beq
\lb e_{ij}^\si, e_{kl}^\theta \rb = \de_{jk} e_{il}^{\si+\theta} -
(-1)^{\si\theta} \de_{il} e_{kj}^{\si+\theta},
\eeq
where $\si,\theta\in\Z_2=\{\bar 0, \bar 1\}$, and $i,j,k,l\in\{0,1, 
\ldots,n\}$. In this paper, we shall use $\lb\,,\,\rb$ for the
Lie superalgebra bracket, and write explicitly $[\,,\,]$ (
resp.\ $\{\,,\,\}$) if this
stands for a commutator (resp.\ anti-commutator).
Let $E_{ij}$ denote the $(n+1)\times(n+1)$ matrix with 1 in position
$(i,j)$ and 0 elsewhere (indices running from 0 to $n$), then the
defining representation of $q(n+1)$ is given by the map
\beq
e_{ij}^{\bar 0} \rightarrow 
\left( \begin{array}{cc} E_{ij}&0\\0&E_{ij} \end{array} \right),\qquad
e_{ij}^{\bar 1} \rightarrow 
\left( \begin{array}{cc} 0&E_{ij}\\E_{ij}&0 \end{array} \right)
\eeq
onto matrices of order $2(n+1)$. 

Following the definition of Cartan subalgebra as a maximal nilpotent
subalgebra coinciding with its own normalizer, 
a Cartan subalgebra $H$ of $G=q(n+1)$ is given by 
$H=H_{\bar 0} \oplus H_{\bar 1}$, where
$H_{\bar 0}={\rm span}\{e_{ii}^{\bar 0}\;|\; i=0,1,\ldots n\}$ and
$H_{\bar 1}={\rm span}\{e_{ii}^{\bar 1}\;|\; i=0,1,\ldots n\}$
\cite{penkov0,dict}. 
This subalgebra is not abelian, and since the elements of $H_{\bar 1}$
are odd, the root generators $e_{ij}^{\bar 1}$ are not eigenvectors
of $H_{\bar 1}$. Therefore, to give a root decomposition of $G$ it
is more convenient~\cite{dict} 
to work with the abelian subalgebra $H_{\bar 0}$.
The dual space $H_{\bar 0}^*$
has the standard basis $\{ \ep_0, \ep_1, \ldots, \ep_n\}$, in terms
of which the roots of $G$ can be described.
The elements $e_{ij}^{\bar 0}$ ($i\ne j$) are even root vectors 
corresponding to the root $\ep_i-\ep_j$;
the elements $e_{ij}^{\bar 1}$ ($i\ne j$) are odd root vectors 
also corresponding to the root $\ep_i-\ep_j$;
the elements $e_{ii}^{\bar 1}$ from $H_{\bar 1}$
can then be interpreted as odd root
vectors corresponding to the root 0. 
Note that every root $\ep_i-\ep_j$ ($i\ne j$) has multiplicity~2 
(counting once as even and once as odd root). 
Let, as usual, $\De$ be the set of all roots, $\De^{\bar 0}$ (resp.\
$\De^{\bar 1}$) be the set of even (resp.\ odd) roots~:
\beq
\De^{\bar 0} = \{ \ep_i -\ep_j \;|\; 0\leq i\ne j \leq n \},\qquad
\De^{\bar 1} = \De^{\bar 0} \cup \{0\}.
\eeq
The positive roots are
\beq
\De_+=\De^{\bar 0}_+ = \De^{\bar 1}_+ = 
\{ \ep_i -\ep_j \;|\; 0\leq i<j \leq n \}.
\eeq
With this choice of positive roots the weights 
$\la=\sum_{i=0}^n \la_i \ep_i \in H^*_{\bar 0}$ are partially ordered by
$\la\leq \mu$ iff $\mu-\la = \sum k_\al \al$ where $\al\in\De_+$
and $k_\al$ are nonnegative integers.
The adjoint representation has $\ep_0-\ep_n$ as highest weight,
with multiplicity~2. The defining representation has highest 
weight $\ep_0$, also with multiplicity~2. 

Let $V$ be a linear space over $\C$, and denote by $gl(V)$
the space of endomorphisms of $V$. A representation $\rh$ is a
linear mapping from $G$ to $gl(V)$ such that
\beq
\rh(\lb x,y \rb ) = \rh(x)\rh(y)- (-1)^{\si\theta} \rh(y)\rh(x),
\qquad \forall x\in G_\si, y\in G_\theta;\ \si,\theta\in\Z_2.
\eeq
Then $V$ is a $G$-module with $xv=\rh(x)v$ for $x\in G$ and $v\in V$.
If, moreover, $V$ is a $\Z_2$-graded linear space, i.e.\
$V=V_{\bar 0} \oplus V_{\bar 1}$, then also $gl(V)$ is naturally
graded, $gl(V)=gl(V)_{\bar 0} \oplus gl(V)_{\bar 1}$, and then
$\rh$ is a graded representation (and $V$ a graded $G$-module)
if $\rh(x)\in gl(V)_\si$ for $x\in G_\si$.
For Lie superalgebras, one usually considers only the graded
modules when studying representation 
theory~\cite{Kac1,Scheunert}. Here, we shall see
that $q(n+1)$ has a class of interesting non-graded modules.

Graded modules of $q(n+1)$ were considered by 
Penkov and Serganova~\cite{penkov0,penkov1,penkov2}.
In particular, they showed that the finite-dimensional irreducible
representations 
$V$ of $q(n+1)$ are characterized by a highest weight
$\la=\sum_i\la_i\ep_i$, such that $\la_i-\la_{i+1}$ is a nonnegative
integer and $\la_i=\la_{i+1}$ implies $\la_i=\la_{i+1}=0$.
The dimension of the highest weight space $V_\la$ ($\la\ne 0$)
is given by~\cite[page~150]{penkov1}
\beq
\dim(V_\la) = 2^{1+[(\#\la-1)/2]},
\eeq
where $\#\la$ is the number of nonzero coordinates $\la_i$, 
and $[t]$ is the integer part of $t$.
For example, for the defining representations with $\la=(1,0,\ldots,0)$
and the adjoint representation with $\la=(1,0,\ldots,0,-1)$, the
highest weight space has dimension~2. 

{}From the physical point of view, it is 
unusual to have a
highest (or lowest) weight with multiplicity greater than~1,
since this is normally
associated to a ``unique vacuum''.

In the present paper we shall show that $q(n+1)$ has a class
of interesting non-graded representations, with a unique highest
weight vector (i.e.\ highest weight multiplicity~1), and which
are also finite-dimensional. Moreover, these representations
can be interpreted as Hermitian Fock representations generated 
by $n$ even and $n$ odd creation operators.

\section{Creation and annihilation operators for $q(n+1)$}

Let $a_i^\xi(\si)$ be the following elements of $q(n+1)$~:
\beq
a_i^+(\si) = e_{i,0}^\si, \qquad
a_i^-(\si) = e_{0,i}^\si, \qquad \si\in\Z_2,\quad i\in\{1,\ldots,n\}.
\label{Q0}
\eeq
It is easy to verify that these operators satisfy the
following relations~:
\bea
&& \lb a_i^-(\si), a_j^-(\theta) \rb =
\lb a_i^+(\si), a_j^+(\theta) \rb = 0, \label{Q1}\\
&& \lb \lb a_i^+(\si), a_j^-(\theta) \rb, a_k^+(\omega) \rb =\nn\\
&&\qquad\qquad \de_{jk} a^+_i(\si+\theta+\omega) + (-1)^{\si\theta+\theta\omega+\omega\si}
\de_{ij} a^+_k(\si+\theta+\omega);\label{Q2}\\
&& \lb \lb a_i^+(\si), a_j^-(\theta) \rb, a_k^-(\omega) \rb =\nn\\
&&\qquad\qquad -(-1)^{\si\theta} \de_{ij} a^-_k(\si+\theta+\omega) 
 - (-1)^{\theta\omega+\omega\si}
\de_{ik} a^-_j(\si+\theta+\omega),\label{Q3}
\eea
where $\si,\theta,\omega\in\Z_2$, $i,j,k\in\{1,\ldots,n\}$.
It is convenient to introduce the following notational difference
between the even and odd operators~:
\beq
b_i^{\pm} = a_i^{\pm}(\bar 0), \qquad
f_i^{\pm} = a_i^{\pm}(\bar 1).
\eeq
The operators $b_i^+,\ f_i^+$ (resp.\  $b_i^-,\ f_i^-$) shall be
referred to as creation (resp.\ annihilation) operators for the Lie
superalgebra $q(n+1)$ (even though they generate only the subalgebra
$sq(n+1)$).

A definition of creation and annihilation operators (or
generators) of a simple Lie (super)algebra ${\cal L}$ and of
the related Fock spaces was given in~\cite[\S 2]{palevSA}. The
motivation for introducing such operators stems from the
observation that the creation and annihilation operators (CAO's)
of certain algebras have a direct physical significance. We have
in mind the para-Fermi and the para-Bose operators, which
generalize the statistics of spinor and tensor fields in quantum
field theory~\cite{Green}. Any $n$ pairs of parafermions are CAO's
of the orthogonal Lie algebra 
$so(2n+1)\equiv B_n$~\cite{Kamefuchi,Ryan}. 
The parabosons do not generate a Lie
algebra, they generate a Lie superalgebra~\cite{Omote}. Any
$n$ pairs of them are CAO's of the orthosymplectic Lie
superalgebra $osp(1,2n)\equiv B(0,n)$~\cite{Ganchev}.

In~\cite{palevD} the question was raised whether each
simple Lie (super)algebra can be generated by creation and
annihilation generators. The answer is positive for all algebras
{}from the classes $A$, $B$, $C$ and $D$ of simple Lie 
algebras~\cite{palevC} and for some Lie superalgebras. 
So far however only the CAO's
and the Fock representations of $sl(n+1)$ 
($A$-statistics)~\cite{palevC,palevA} and of Lie superalgebras 
$sl(1/n)$ ($A-$superstatistics)~\cite{palevSA} were studied in 
somewhat greater detail. The present paper is another contribution
along this line for the Lie superalgebra $q(n+1)$.

Coming back to the CAO's (\ref{Q0}) we note that the operators
$b_i^{\pm}$ satisfy the
relations of $A$-statistics~\cite{palevA}, whereas the operators
$f_i^{\pm}$ satisfy the relations of $A$-superstatistics~\cite{palevSA}.
Here, we shall refer to the combined relations (\ref{Q1})-(\ref{Q3}) as 
$Q$-statistics.
Clearly, the linear envelope of 
\beq
\{ a_i^\xi(\si), \lb a_i^\xi(\si),a_j^\eta(\theta) \rb |
\xi,\eta\in\{+,-\},\, \si,\theta\in\Z_2,\, i,j\in\{1,\ldots,n\} \}
\eeq
is equal to the Lie superalgebra $sq(n+1)$.

\section{Fock space for $q(n+1)$}

We shall define a Fock space for $G=q(n+1)$ using an induced module.
First of all, from the commutation relations of $q(n+1)$ it is
straightforward to see that $q(1)={\rm span} \{ e_{00}^{\bar 0},
e_{00}^{\bar 1} \}$ and $q(n)={\rm span} \{ e_{i,j}^{\si} \;|\;
i,j=1,2,\ldots,n;\ \si=\bar 0,\bar 1 \}$ 
are subalgebras of $q(n+1)$ with 
$\lb q(1), q(n) \rb =0$. So consider the subalgebra
\beq
\tilde G = q(1) \oplus q(n).
\eeq
Let 
\bea
P&=&{\rm span} \{b_1^-,\ldots,b_n^-, f_1^-,\ldots,f_n^-\},\\ \nn
N&=&{\rm span} \{b_1^+,\ldots,b_n^+, f_1^+,\ldots,f_n^+\};
\eea
these are two abelian subalgebras of $G$.
Then $G=\tilde G +P+N$, where the sum is direct as linear spaces.
Since $\lb \tilde G,P\rb = P$, $\tilde G+P$ is also a subalgebra
of $G$.

The Lie superalgebra $q(1)$ has basis elements $e_{00}^{\bar 0}$
and $e_{00}^{\bar 1}$ with supercommutation relations
\beq
[e_{00}^{\bar 0},e_{00}^{\bar 0}]=0, \qquad
[e_{00}^{\bar 0},e_{00}^{\bar 1}]=0, \qquad
\{e_{00}^{\bar 1},e_{00}^{\bar 1} \}=2e_{00}^{\bar 0}.
\eeq
Clearly, this Lie superalgebra has one-dimensional irreducible 
modules $\C v_0$ characterized by a number $p$, with action
\beq
e_{00}^{\bar 0} v_0=p\; v_0,\qquad e_{00}^{\bar 1} v_0=\sqrt{p}\; v_0.
\label{v0}
\eeq
In principle $p$ can be any complex number, but later we shall be
interested only in the case that $p$ is a positive real number.
The $q(1)$-module $\C v_0$ can be extended to a $\tilde G= q(1) \oplus
q(n)$-module by letting $x v_0=0$ for all $x\in q(n)$.
Requiring that $x v_0=0$ for every $x\in P$ it becomes a
$(\tilde G+P)$-module.

We now define the following induced $G$-module~:
\beq
\bar V_p = {\rm Ind}_{\tilde G +P}^G \C v_0 \cong 
U(G) \otimes_{\tilde G+P} \C v_0.
\eeq
By the Poincar\'e-Birkhoff-Witt theorem for Lie superalgebras, we have
\beq
\bar V_p \cong U(N) \otimes \C v_0.
\eeq
Thus a basis of $\bar V_p$ is given by the elements
\bea
&&|p;k_1,l_1,k_2,l_2,\ldots,k_n,l_n\rangle =
(b_1^+)^{k_1}(f_1^+)^{l_1}(b_2^+)^{k_2}(f_2^+)^{l_2}\cdots 
(b_n^+)^{k_n}(f_n^+)^{l_n} v_0,\nn\\
&&\qquad l_i\in\{0,1\},\; k_i=0,1,2,\ldots .
\label{pkl}
\eea
What are the identities that hold in this representation space 
$\bar V_p$? First of all, note that $\bar V_p$ has a unique highest 
weight equal to $p\ep_0$, corresponding to the unique (up to
a factor) highest weight vector $v_0$. So the highest weight
has multiplicity~1. On the other hand, (\ref{v0}) shows that $v_0$
is not an even nor an odd vector, i.e.\ the $G$-module $\bar V_p$
is not graded. Note that the weight of (\ref{pkl}) is given by
\beq
p\ep_0 + \sum_{i=1}^n (k_i+l_i)(\ep_i-\ep_0).
\eeq
Secondly, the vector $v_0$ can genuinely be called a ``vacuum vector''
since it satisfies
\beq
b_i^- v_0 = f_i^- v_0 =0, \qquad (i=1,2,\ldots,n).
\label{bfm}
\eeq
Furthermore, the following relations are valid~:
\bea
&& b_i^- b_j^+ v_0 = \de_{ij} p\; v_0, \qquad
f_i^- f_j^+ v_0 = \de_{ij} p\; v_0,\nn \\ 
&& f_i^- b_j^+ v_0 = \de_{ij} \sqrt{p}\; v_0, \qquad
b_i^- f_j^+ v_0 = \de_{ij} \sqrt{p}\; v_0. \label{minplus}
\eea
Note that relations (\ref{bfm}) and (\ref{minplus}) 
are also sufficient to define
the representation $\bar V_p$.

In order to call $\bar V_p$ a Fock space, one further condition
should be satisfied, namely it should be a Hilbert space consistent
with the adjoint operation~\cite{palevA,palevSA}
\beq
(b_i^\pm)^\dagger = b_i^\mp, \qquad 
(f_i^\pm)^\dagger = f_i^\mp.
\eeq
Thus we define a Hermitian form on $\bar V_p$ by
\beq
\langle v_0|v_0 \rangle =1,\;
\langle b_i^+v|w \rangle=\langle v|b_i^-w \rangle,\;
\langle f_i^+v|w \rangle=\langle v|f_i^-w \rangle,\qquad
v,w\in \bar V_p.
\label{form}
\eeq
In general $\bar V_p$ is not a Hilbert space. However,
we shall see that if $p$ is a positive integer, $\bar V_p$ has
a quotient space which is a Hilbert space.
Indeed, if $p$ is a positive integer, the space $\bar V_p$ is shown to
have a maximal submodule $M_p$. Then the quotient module 
$V_p=\bar V_p / M_p$ is an irreducible $G$-module.
The Hermitian form is zero on $M_p$ and on $V_p$ it induces a 
positive definite metric. Thus $V_p$ can genuinely be called a
Fock space representation of $q(n+1)$.

The Lie superalgebra $q(n+1)$ contains a one-dimensional center,
\beq
I=\sum_{i=0}^n e_{ii}^{\bar 0}.
\eeq
Since $Iv_0=pv_0$, and $\bar V_p$ (or $V_p$) is generated by $v_0$, 
it follows that $Iv=pv$ for every $v$ in $\bar V_p$ (or $V_p$).

As we shall see, the structure of $V_p$ or $\bar V_p$ is far from trivial.
Before turning to the general case, let us first consider the
low rank case of the Lie superalgebra $G=q(2)$.

\section{Fock space for $q(2)$}

Since $n=1$ there is only one index for the creation and annihilation
operators, so we shall simply denote $b^\pm_1, f^\pm_1$ by $b^\pm, f^\pm$.

The representation space $\bar V_p$ is spanned by the following vectors
(notation of (\ref{pkl}))~:
\bea
v_k= |p;k,0\rangle = (b^+)^k v_0,\; k=0,1,\ldots;\\ \nn
w_k= |p;k-1,1\rangle = (b^+)^{k-1} f^+ v_0,\; k=1,2,\ldots.
\eea
The following actions of the annihilation operators on $v_k$ and $w_k$
can be computed using the triple relations (\ref{Q1}) and 
(\ref{Q2}), (\ref{Q3})~:
\bea
&& b^- v_k = k(p-k+1) v_{k-1},\label{q21}\\
&& f^- v_k = k \sqrt{p}\; v_{k-1} -k(k-1) w_{k-1},\label{q22}\\
&& b^- w_k = \sqrt{p}\; v_{k-1} + (k-1)(p-k) w_{k-1},\label{q23}\\
&& f^- w_k = p v_{k-1} - (k-1)\sqrt{p}\; w_{k-1}. \label{q24}
\eea
It is not difficult to verify the following, using the earlier defined
metric on $\bar V_p$ and (\ref{q21})~:
\beq
\langle v_k | v_k \rangle = k! p(p-1)\cdots (p-k+1).
\eeq
For fixed $p$, this expression can take positive and negative values,
depending upon the value of $k$. Hence, $\bar V_p$ itself is not
a Hilbert space representation. 
Next, let us investigate whether $\bar V_p$ is irreducible. Using
the relations (\ref{form}), (\ref{q21})-(\ref{q24}), 
and induction, one can show that
\bea
&& \langle v_k | v_l \rangle = \de_{kl} k! p(p-1)\cdots (p-k+1),
\label{q2f1}\\
&& \langle w_k | w_l \rangle = \de_{kl} (k-1)! p(p-1)\cdots (p-k+1),
\label{q2f2}\\
&& \langle v_k | w_l \rangle = \de_{kl} k! p(p-1)\cdots (p-k+1)/\sqrt{p}.
\label{q2f3}
\eea
Now (\ref{q2f1}) and (\ref{q2f2}) imply that
\bea
&& (b^-)^k v_k = k! p(p-1)\cdots (p-k+1)\, v_0,\\ \nn
&& (b^-)^{k-1}f^- w_k = (k-1)! p(p-1)\cdots (p-k+1)\, v_0.
\eea
Thus if $p$ is not a positive integer these coefficients are not
zero, implying that the vectors $v_k$ and $w_k$ cannot belong to
a submodule of $\bar V_p$ (apart from the trivial submodule $\bar V_p$).
In other words, if $p$ is not a positive integer, $\bar V_p$ 
is irreducible. 

Let us now consider the interesting case that $p$ is a positive integer.
Then $\bar V_p$ has a maximal submodule $M_p$. Since $\bar V_p$ is
a module generated by a highest weight vector, 
the submodule is generated by primitive weight vectors, 
so let us determine when a weight vector $v_k+\be w_k$
is primitive. The conditions $b^-(v_k+\be w_k)=0$ and
$f^-(v_k+\be w_k)=0$ lead to one solution only, namely $k=p$ and 
$\be=-\sqrt{p}$. Thus $v_p-\sqrt{p}\; w_p$ is a primitive vector generating
the submodule $M_p$. The quotient module $V_p=\bar V_p/M_p$ is therefore
a finite-dimensional module. A set of basis vectors of $V_p$,
together with the corresponding weight, is given by
\beq
\begin{array}{lll}
v_0 && p\ep_0\\
v_1, w_1 && (p-1)\ep_0+\ep_1\\
v_2, w_2 && (p-2)\ep_0+2\ep_1\\
\vdots && \vdots \\
v_{p-1}, w_{p-1} && \ep_0+(p-1)\ep_1\\
v_p+\sqrt{p}\; w_p && p\ep_1.
\end{array}
\eeq
The top and bottom weight appear with multiplicity~1, the other weights
have multiplicity~2. From the weight structure one can determine the
decomposition of this finite-dimensional $q(2)$ module with respect
to the subalgebra $gl(2)\subset q(2)$~:
\beq
V_p \rightarrow (p,0) \oplus (p-1,1),\qquad\qquad (p>1).
\eeq
So $V_p$ splits into two irreducible $gl(2)$ modules, both of which have
been labelled by their highest weight (in the $(\ep_0,\ep_1)$-basis). 
So for $p>1$ these $q(2)$ representations could be referred to
as ``dispin'' representations, similar to the known dispin representations
of $osp(1,2)$~\cite{SNR,Hughes}.
For $p=1$, the decomposition
is simply $V_p \rightarrow (p,0)$. The dimension follows easily~:
\beq
\dim V_p = 2p.
\eeq

It is possible to give an orthonormal basis for $V_p$, in terms of
the above basis vectors $v_k$, $w_k$. Since
\bea
&&\langle v_k\pm \sqrt{k} w_k | v_k\pm \sqrt{k} w_k \rangle =
2(1\pm\sqrt{k/p}) k!p!/(p-k)!,\\
&&\langle v_k+ \sqrt{k} w_k | v_k- \sqrt{k} w_k \rangle =0,
\eea
one can define
\bea
&&\phi_k = \left( (p-k)!\over 2k!p!(1+\sqrt{k/p}) \right)^{1/2} 
(v_k+ \sqrt{k} w_k),\;(k=1,\ldots,p)\\
&&\psi_k = \left( (p-k)!\over 2k!p!(1-\sqrt{k/p}) \right)^{1/2} 
(v_k- \sqrt{k} w_k),\;(k=0,\ldots,p-1).
\eea
These vectors are orthonormal~:
\beq
\langle \phi_k | \phi_l \rangle = \langle \psi_k | \psi_l \rangle = 
\de_{kl},\qquad \langle \phi_k | \psi_l \rangle = 0.
\eeq
The action of the creation and annihilation operators on this basis
can be computed. We have~:
\bea
f^+ \phi_k &=& {1\over 2}\left[ (\sqrt{p}-\sqrt{k})(\sqrt{p}+\sqrt{k+1})
\right]^{1/2} \phi_{k+1} \nn\\
&&-{1\over 2}\left[ (\sqrt{p}-\sqrt{k})(\sqrt{p}-\sqrt{k+1})
\right]^{1/2} \psi_{k+1}\\
f^+ \psi_k &=& {1\over 2}\left[ (\sqrt{p}+\sqrt{k})(\sqrt{p}+\sqrt{k+1})
\right]^{1/2} \phi_{k+1} \nn\\
&&-{1\over 2}\left[ (\sqrt{p}+\sqrt{k})(\sqrt{p}-\sqrt{k+1})
\right]^{1/2} \psi_{k+1}\\
b^+ \phi_k &=& {1\over 2}(\sqrt{k+1}+\sqrt{k})
\left[ (\sqrt{p}-\sqrt{k})(\sqrt{p}+\sqrt{k+1})
\right]^{1/2} \phi_{k+1} \nn\\
&&+{1\over 2}(\sqrt{k+1}-\sqrt{k})
\left[ (\sqrt{p}-\sqrt{k})(\sqrt{p}-\sqrt{k+1})
\right]^{1/2} \psi_{k+1}\\
b^+ \psi_k &=& {1\over 2}(\sqrt{k+1}-\sqrt{k})
\left[ (\sqrt{p}+\sqrt{k})(\sqrt{p}+\sqrt{k+1})
\right]^{1/2} \phi_{k+1} \nn\\
&&+{1\over 2}(\sqrt{k+1}+\sqrt{k})
\left[ (\sqrt{p}+\sqrt{k})(\sqrt{p}-\sqrt{k+1})
\right]^{1/2} \psi_{k+1}.
\eea
The action of the annihilation operators follows immediately from 
$b^-=(b^+)^\dagger$, $f^-=(f^+)^\dagger$. For example,
\beq
f^- \phi_k = {1\over 2}\left[ (\sqrt{p}-\sqrt{k-1})(\sqrt{p}+\sqrt{k})
\right]^{1/2} \phi_{k-1} +
{1\over 2}\left[ (\sqrt{p}+\sqrt{k-1})(\sqrt{p}+\sqrt{k})
\right]^{1/2} \psi_{k-1}.
\eeq

\section{Structure of the module $\bar V_p$}

In this section we return to the general case $q(n+1)$. By 
calculating the action of creation and annihilation operators on
basis vectors of the induced module $\bar V_p$, the way is
prepared to determine the structure of the irreducible quotient
module $V_p$ of $\bar V_p$.

Recall that a basis for $\bar V_p$ is given by the vectors
\bea
|p;{\pmb k}, {\pmb l } \rangle
&=&|p;k_1,l_1,k_2,l_2,\ldots,k_n,l_n\rangle =
(b_1^+)^{k_1}(f_1^+)^{l_1}(b_2^+)^{k_2}(f_2^+)^{l_2}\cdots 
(b_n^+)^{k_n}(f_n^+)^{l_n} v_0,\nn\\
&&\qquad l_i\in\{0,1\},\; k_i=0,1,2,\ldots . \label{pkkll}
\eea
Since all creation operators $b_i^+$ and $f_i^+$ supercommute,
a different order of the creation operators in (\ref{pkkll}) can produce
only a sign change.

In the standard basis, the weight of the vector 
$|p;{\pmb k}, {\pmb l } \rangle$ is given by
\beq
{\rm weight}(|p;{\pmb k}, {\pmb l } \rangle) =
p\ep_0+\sum_{i=1}^n (k_i+l_i)(\ep_i-\ep_0)=
\left( p-\sum_{i=1}^n(k_i+l_i), k_1+l_1,\ldots, k_n+l_n\right).
\eeq
Thus every weight of $\bar V_p$ is of the form
\beq
\la_m = (p-\sum_{i=1}^n m_i, m_1,\ldots, m_n),\qquad m_i=0,1,2,\ldots ;
\label{la_m}
\eeq
conversely, every weight of the form (\ref{la_m}) 
is a weight of $\bar V_p$.
Since the $l_i$ in (\ref{pkkll}) are either 0 or 1, it follows that
the multiplicity of the weight $\la_m$ in $\bar V_p$ is given by
\beq
\mult_{\bar V_p} (p-\sum_{i=1}^n m_i, m_1,\ldots, m_n)=
2^{\ga(m_1)+\cdots +\ga(m_n)},
\eeq
where $\ga(m_i)=0$ if $m_i=0$ and $\ga(m_i)=1$ if $m_i\ne 0$.

The action of the creation operators on the basis of ${\bar V_p}$
is very simple. It will be convenient to denote the basis vectors
in the right hand side of such actions only by means of the labels
that are effectively changed by the action. So, instead of
writing
\beas
&& b^+_j |p;{\pmb k}, {\pmb l } \rangle = 
|p;k_1,l_1,\ldots,k_j+1,l_j,\ldots,k_n,l_n\rangle,\\
&& f^+_j |p;{\pmb k}, {\pmb l } \rangle = \de_{l_j,0}
(-1)^{l_1+\cdots+l_{j-1}}
|p;k_1,l_1,\ldots,k_j,l_j+1,\ldots,k_n,l_n\rangle,
\eeas
we abbreviate this to~:
\bea
&& b^+_j |p;{\pmb k}, {\pmb l } \rangle = |k_j+1 \rangle,\label{bj+}\\
&& f^+_j |p;{\pmb k}, {\pmb l } \rangle = \de_{l_j,0}
(-1)^{l_1+\cdots+l_{j-1}} |l_j+1 \rangle. \label{fj+}
\eea
The action of the annihilation operators is more complicated, and here
the notational convention just introduced will be very useful.

\begin{prop}
The action of the annihilation operators in the module $\bar V_p$
is given by~:
\bea
f_j^- |p;{\pmb k}, {\pmb l } \rangle &=&
(-1)^{l_1+\cdots+l_{j-1}} l_j \bigl(p+1+k_j-\sum_{i=1}^n(k_i+l_i)\bigr)
 |l_j-1\rangle \nn\\
&& + (-1)^{l_1+\cdots+l_n} k_j \sqrt{p} |k_j-1\rangle 
- (-1)^{l_1+\cdots+l_{j-1}} \de_{l_j,0} k_j(k_j-1) |k_j-2,l_j+1\rangle\nn\\
&& - \sum_{i=1 \atop i\ne j}^n (-1)^{l_1+\cdots+l_{i-1}} \de_{l_i,0}
k_i k_j |k_j-1, k_i-1,l_i+1 \rangle \nn\\
&& + \sum_{i=1 \atop i\ne j}^n (-1)^{l_1+\cdots+l_{i-1}} 
l_i k_j |k_j-1, k_i+1,l_i-1 \rangle; \label{fj-}
\eea
\bea
(-1)^{l_{j+1}+\cdots+l_n}b_j^- |p;{\pmb k}, {\pmb l } \rangle &=&
(-1)^{l_{j+1}+\cdots+l_n} k_j \bigl(p+1-l_j-\sum_{i=1}^n(k_i+l_i)\bigr)
|l_j-1\rangle \nn\\
&&+ l_j \sqrt{p} |l_j-1\rangle 
+ \sum_{i=1 \atop i\ne j}^n (-1)^{l_i+\cdots+l_n}\th_{ij}\de_{l_i,0}
k_i l_j |l_j-1, k_i-1,l_i+1 \rangle \nn\\
&& - \sum_{i=1 \atop i\ne j}^n (-1)^{l_i+\cdots+l_n} \th_{ij}
l_i l_j |l_j-1, k_i+1,l_i-1 \rangle, \label{bj-}
\eea
where $\th_{ij}=+1$ if $i<j$ and $\th_{ij}=-1$ if $i>j$.
\end{prop}

\noindent {\bf Proof.} We shall sketch the proof for the action of
$f_j^-$; that for $b_j^-$ is similar. The proof uses induction on $n$.

As a first step, the action of $f_1^+$ will be determined.
Denote $y_i=(b_i^+)^{k_i}(f_i^+)^{l_i}$. Then,
\beq
f_1^- |p;{\pmb k}, {\pmb l } \rangle = 
f_1^- y_1 y_2\cdots y_n v_0 = \lb f_1^-,y_1\rb y_2\cdots y_n v_0 +
(-1)^{l_1} y_1 f_1^- y_2 \cdots y_n v_0. \label{54}
\eeq
Since the weight of $f_1^- y_2 \cdots y_n v_0$, which is
$\bigl(p+1-\sum_{i=2}^n(k_i+l_i),-1,k_2+l_2,\ldots, k_n+l_n\bigr)$, is
not of the form (\ref{la_m}), the vector cannot belong to $\bar V_p$, so 
the second term in (\ref{54}) has to be zero. 
Using $\lb f_1^-,b_1^+\rb=e_{00}^{\bar 1}- e_{11}^{\bar 1}$ and
$\lb f_1^-,f_1^+\rb=e_{00}^{\bar 0}+ e_{11}^{\bar 0}$, one finds
\beq
\lb f_1,y_1\rb = \sum_{r=1}^{k_1-1} (b_1^+)^r 
(e_{00}^{\bar 1}- e_{11}^{\bar 1}) (b_1^+)^{k_1-r-1}(f_1^+)^{l_1} +
l_1 (b_1^+)^{k_1}(e_{00}^{\bar 0}+ e_{11}^{\bar 0}).
\eeq
{}From the weight of $y_i$ and $v_0$ one obtains
\beq
(e_{00}^{\bar 0}+ e_{11}^{\bar 0})y_2\cdots y_n v_0 =
(p-k_2-l_2-\cdots-k_n-l_n) y_2\cdots y_n v_0.
\eeq
Next we need to determine the action of $e_{00}^{\bar 1}$ and 
$e_{11}^{\bar 1}$ on vectors of the form $y'_1 y_2\cdots y_n v_0$,
where $y'_1=(b_1^+)^{k'_1}(f_1^+)^{l_1}$ with $k'_1=k_1-r-1$~:
\beas
e_{00}^{\bar 1} y'_1 y_2\cdots y_n v_0 &=& 
\lb e_{00}^{\bar 1}, y'_1 y_2\cdots y_n\rb v_0 +
(-1)^{l_1+\cdots l_n} y'_1 y_2\cdots y_n  \sqrt{p}\; v_0\\ \nn
&=& \lb e_{00}^{\bar 1}, y'_1\rb y_2\cdots y_n v_0 +
\sum_{i=2}^n y'_1\cdots y_{i-1}\lb e_{00}^{\bar 1}, y_i\rb
  y_{i+1}\cdots y_n v_0  \\ \nn
&& \qquad\qquad + (-1)^{l_1+\cdots l_n} \sqrt{p}\; y'_1 y_2\cdots y_n v_0.
\eeas
Every term in this expression can be determined explicitly using
\beq
\lb e_{00}^{\bar 1}, y_i\rb = -\de_{l_i,0} k_i (b_i^+)^{k_i-1} f_i^+ +
l_i (b_i^+)^{k_i+1},
\eeq
which follows from the supercommutator of $e_{00}^{\bar 1}$
with $b_i^+$ and $f_i^+$. 
Similarly, for the action of $e_{11}^{\bar 1}$ one finds
\beq
e_{11}^{\bar 1} y'_1 y_2\cdots y_n v_0 = 
\lb e_{11}^{\bar 1}, y'_1 \rb y_2\cdots y_n v_0 ,
\eeq
since $\lb e_{11}^{\bar 1}, y_i \rb =0$ for $i>1$ and 
$e_{11}^{\bar 1} v_0=0$. Using the supercommutator of $e_{11}^{\bar 1}$
with $b_1^+$ and $f_1^+$, there comes
\beq
\lb e_{11}^{\bar 1}, y'_1 \rb = \de_{l_1,0} k'_1 (b_1^+)^{k'_1-1}
f_1^+ + l_1 (b_1^+)^{k'_1+1}.
\eeq
Collecting now all contributions yields the action of $f_1^+$ on
$|p;{\pmb k}, {\pmb l } \rangle$, as given in the proposition.

In the second step, we use induction in the following way. 
First observe that for $j=2$,
\beq
f_2^- y_1 y_2\cdots y_n v_0 = \lb f_2^-, y_1\rb y_2\cdots y_n v_0
+ (-1)^{l_1} y_1 f_2^- y_2\cdots y_n v_0.
\eeq
But the action $f_2^- y_2\cdots y_n v_0$ is formally the same as
the (known) action $f_1^- y_1 y_2\cdots y_n v_0$, by relabelling of
indices. More generally, for $j\geq 2$,
\beq
f_j^- y_1 y_2\cdots y_n v_0 = \lb f_j^-, y_1\rb y_2\cdots y_n v_0
+ (-1)^{l_1} y_1 f_j^- y_2\cdots y_n v_0.\label{61}
\eeq
Once again, $f_j^- y_2\cdots y_n v_0$ can formally be reduced to
$f_{j-1}^- y_1\cdots y_{n-1} v_0$, on which one uses the induction
hypothesis. So what remains to be determined is the first term 
in (\ref{61}).
Since $[f_j^-,b_1^+]$ commutes with $b_1^+$ (see (\ref{Q2})), one
finds
\beq
\lb f_j^-,(b_1^+)^{k_1}(f_1^+)^{l_1}\rb =
(-1)^{l_1} k_1 (b_1^+)^{k_1-1}(f_1^+)^{l_1} [ f_j^-,b_1^+ ] +
l_1 (b_1^+)^{k_1} \{ f_j^-, f_1^+ \}.
\eeq
To calculate $[ f_j^-,b_1^+ ] y_2\cdots y_n v_0$, one uses again
the triple relations (\ref{Q2}) and (\ref{Q3})~:
\beq
[ f_j^-,b_1^+ ] y_2\cdots y_n v_0 = (-1)^{l_2+\cdots+l_{j-1}} 
y_2\cdots y_{j-1} \lb [ f_j^-,b_1^+ ], y_j \rb y_{j+1}\cdots y_n v_0;
\eeq
furthermore
\beq
\lb [ f_j^-,b_1^+ ], y_j \rb = -k_j f_1^+ (b_j^+)^{k_j-1} (f_j^+)^{l_j}
-l_j (b_j^+)^{k_j} b_1^+,
\eeq
where again the triple relations have been used. In a similar way,
the action $\{ f_j^-,f_1^+ \} y_2\cdots y_n v_0$ is determined. 
Together they yield the first term in the right hand side of (\ref{61}).
Combining then the coefficients of all identical vectors in the 
right hand side of (\ref{61}) proves the proposition. \mybox

Now we wish to determine the vectors in $\bar V_p$ that are 
annihilated by all annihilation operators, i.e. by
\beq
P={\rm span} \{b_1^-,\ldots,b_n^-, f_1^-,\ldots,f_n^-\}.
\eeq
For this purpose, and inspired by the right hand sides 
in (\ref{bj-}) and (\ref{fj-}), we
introduce the following weight vectors
(still using the notational convention introduced in 
(\ref{bj+})-(\ref{fj+}))~:
\bea
X(p;{\pmb k},{\pmb l}) &=& (-1)^{l_1+\cdots+l_n} \sqrt{p}
|p;{\pmb k}, {\pmb l } \rangle -
\sum_{i=1}^n (-1)^{l_1+\cdots+l_i} l_i |k_i+1,l_i-1 \rangle \nn \\ 
&& - \sum_{i=1}^n (-1)^{l_1+\cdots+l_i} \de_{l_i,0} 
k_i |k_i-1,l_i+1 \rangle .
\eea
Then we have~:

\begin{prop}
For the weight vectors $|p;{\pmb k},{\pmb l} \rangle$ and
$X(p;{\pmb k},{\pmb l})$ the following
equalities hold (again we use the convention
that in the right hand side only the labels $k_i$ and $l_i$ that
are effectively changed are withheld)~:
\bea
b_j^- |p;{\pmb k},{\pmb l} \rangle &=& k_j(p+1-\sum_{i=1}^n(k_i+l_i))
|k_j-1\rangle + (-1)^{l_1+\cdots+l_{j-1}}l_j X(p;l_j-1),\label{bj-p}\\
f_j^- |p;{\pmb k},{\pmb l} \rangle &=& 
(-1)^{l_1+\cdots+l_{j-1}}l_j(p+1-\sum_{i=1}^n(k_i+l_i))
|l_j-1\rangle + k_j X(p;k_j-1),\label{fj-p}\\
b_j^+ X(p;{\pmb k},{\pmb l}) &=& X(p;k_j+1)+(-1)^{l_1+\cdots+l_j}
\de_{l_j,0} |l_j+1\rangle,\label{bj+x}\\
f_j^+ X(p;{\pmb k},{\pmb l}) &=& -\de_{l_j,0} 
(-1)^{l_1+\cdots+l_{j-1}} X(p;l_j+1)+ |k_j+1\rangle,\label{fj+x}\\
b_j^- X(p;{\pmb k},{\pmb l}) &=& k_j(p-\sum_{i=1}^n(k_i+l_i))
X(p;k_j-1),\label{bj-x}\\
f_j^- X(p;{\pmb k},{\pmb l}) &=& (-1)^{l_1+\cdots+l_{j-1}}
l_j(p-\sum_{i=1}^n(k_i+l_i)) X(p;l_j-1). \label{fj-x}
\eea
\end{prop}
The proof is by direct computation, using proposition~1.
The first two relations are just a reformulation of the equalities
in proposition~1. The next two relations follow immediately from
the definition of $X(p;{\pmb k},{\pmb l})$ and the actions 
(\ref{bj+})-(\ref{fj+}).
The last two relations are the most difficult to verify. 
They follow from the action of the annihilation operators on each
part of $X(p;{\pmb k},{\pmb l})$, using proposition~1, and then
collecting the terms according to equal weight vectors 
$|p;{\pmb k},{\pmb l} \rangle$. This computation is long but
straightforward, and will not be given in detail here. \mybox

Note that all vectors $X(p;{\pmb k},{\pmb l})$ with 
$\sum_i (k_i+l_i)=p$ are annihilated by $P$. 

In order to understand the linear (in)dependance of the newly
introduced weight vectors $X(p;{\pmb k},{\pmb l})$, it is
convenient to prove first an interesting lemma about the determinant and
rank of a matrix. Let $r$ be a positive integer, and consider the
$2^r\times 2^r$ matrix $A$, where the rows and columns of $A$
are labelled by the binary sequences ${\pmb l}=(l_1,\ldots,l_r)$
($l_i\in\{0,1\}$) of length $r$ in reverse binary
order. The elements $a_{{\pmb l},{\pmb l'}}$ of $A\equiv 
A(s;t_1,\ldots,t_r)$ are 
as follows
\bea
&& {\rm if}\ {\pmb l}={\pmb l'}\ {\rm then}\ 
a_{{\pmb l},{\pmb l'}} = s,\nn \\ \nn
&& {\rm if}\ {\pmb l}\ {\rm and}\ {\pmb l'}\ 
{\rm differ\ in\ only\ position}\ i \ {\rm then} \\ \label{defA}
&& \qquad\qquad {\rm if}\ l_i=0\ {\rm then}\  
a_{{\pmb l},{\pmb l'}} = -(-1)^{l_{i+1}+\cdots+l_r} t_i,\\ \nn
&& \qquad\qquad {\rm if}\ l_i=1\ {\rm then}\  
a_{{\pmb l},{\pmb l'}} = (-1)^{l_i+\cdots+l_r} t_i,\\ \nn
&& {\rm otherwise}\ a_{{\pmb l},{\pmb l'}} = 0.
\eea
Herein, $s$ and $t_1,\ldots,t_r$ are arbitrary real numbers or
variables.

It is constructive to consider an example. For $r=2$, the binary
sequences labelling the rows and columns of $A$ are $(0,0), (1,0),
(0,1), (1,1)$ (in this order); for $r=3$ the sequences are
$(0,0,0), (1,0,0), (0,1,0), (1,1,0), (0,0,1), (1,0,1), (0,1,1), (1,1,1)$.
The matrices take the following form~:
\beq
A(s;t_1,t_2) = \left(
\begin{array}{cccc}
s & -t_1 & -t_2 & 0 \\
-1 & s & 0 & -t_2 \\
-1 & 0 & s & t_1 \\
0 & -1 & 1 & s
\end{array} \right),
\eeq
\beq
A(s;t_1,t_2,t_3) = \left(
\begin{array}{cccccccc}
s & -t_1 & -t_2 & 0 & -t_3 &0&0&0\\
-1 & s & 0 & -t_2 &0&-t_3&0&0\\
-1 & 0 & s & t_1 &0&0&-t_3&0\\
0 & -1 & 1 & s &0&0&0&-t_3\\
-1&0&0&0&s & t_1 & t_2 & 0 \\
0&-1&0&0&1 & s & 0 & t_2 \\
0&0&-1&0&1 & 0 & s & -t_1 \\
0&0&0&-1&0 & 1 & -1 & s
\end{array} \right)
\eeq

\begin{lemm}
The matrix $A(s;t_1,\ldots,t_r)$ ($r>1$) defined above satisfies~:
\begin{itemize}
\item[(a)] the determinant is given by
\beq
\det A(s;t_1,\ldots,t_r) = (s^2-\sum_{i=1}^r t_i)^{2^{r-1}};
\label{detA}
\eeq
\item[(b)]
if the elements $t_i$ are positive real numbers such that 
$\sum_{i=1}^r t_i = s^2$ then $\rank(A(s;t_1,\ldots,t_r))=2^{r-1}$.
\item[(c)] 
\beq
A(s;t_1,\ldots,t_r)\cdot A(-s;t_1,\ldots,t_r)=
(\sum_{i=1}^r t_i-s^2)I,
\label{invA}
\eeq
where $I$ is the identity matrix of order $2^r$.
\end{itemize}
\end{lemm}

\noindent {\bf Proof.}
The proof of (a) is by induction on $r$. Clearly it holds for $r=2$.
By definition, the matrix $A(s;t_1,\ldots,t_r)$ 
can be written in block form as
\beq
A(s;t_1,\ldots,t_r) = \left(
\begin{array}{cc}
A(s;t_1,\ldots,t_{r-1}) & -t_r I \\
-I & -A(-s;t_1,\ldots,t_{r-1})
\end{array} \right),
\label{block}
\eeq
where $I$ is the identity matrix of order $2^{r-1}$. So
by induction we have $\det(A(s;t_1,\ldots,t_{r-1}))=
(s^2-\sum_{i=1}^{r-1} t_i)^{2^{r-2}}$
and $\det(-A(-s;t_1,\ldots,t_{r-1}))=
((-s)^2-\sum_{i=1}^{r-1} t_i)^{2^{r-2}}$.
Then it follows from (\ref{block}) that
\bea
&&\det(A(s;t_1,\ldots,t_{r-1},0))=\nn\\
&&\det(A(s;t_1,\ldots,t_{r-1}))
\det(-A(-s;t_1,\ldots,t_{r-1})) = 
(s^2-\sum_{i=1}^{r-1} t_i)^{2^{r-1}}.
\label{78}
\eea
On the other hand, exchanging $t_i$ and $t_j$ in 
$A(s;t_1,\ldots,t_r)$ corresponds
to a permutation of the rows and corresponding columns of 
$A(s;t_1,\ldots,t_r)$;
a closer examination shows that the signature of such a permutation
is positive. Thus $A(s;t_1,\ldots,t_r)$ 
is invariant for transpositions
of the form $t_i \leftrightarrow t_j$. Therefore, 
$\det A(s;t_1,\ldots,t_r)$
is a symmetric polynomial in the elements $t_i$.
Since the power sum symmetric functions form a basis of the ring
of symmetric polynomials~\cite[p.~24]{Macdonald}, it follows that
\beq
\det A(s;t_1,\ldots,t_{r-1},t_r)= \sum_{k=0}^N \sum_{\kappa\vdash k}
c_\kappa(s) p_\kappa(t_1,\ldots,t_{r-1},t_r),
\label{refnew1}
\eeq
where $N=2^{r-1}$, $\kappa$ is summed over all partitions of $k$,
and $p_\kappa$ is the multiplicative power sum 
function~\cite[p.~24]{Macdonald}.
Combining (\ref{refnew1}) with (\ref{78}) gives
\beq
\det A(s;t_1,\ldots,t_{r-1},0)= \sum_{k=0}^N \sum_{\kappa\vdash k}
c_\kappa(s) p_\kappa(t_1,\ldots,t_{r-1})
= (s^2-p_1(t_1,\ldots,t_{r-1}))^N,
\label{refnew2}
\eeq
since $p_1(t_1,\ldots,t_{r-1})=\sum_{i=1}^{r-1} t_i$. 
Thanks to the linear independence of the $p_\kappa$, the expansion of
the factor to the $N$th power in (\ref{refnew2}) fixes all the
coefficients $c_\kappa(s)$. Substituting these back into (\ref{refnew1}),
it follows that we must have
\beq
\det A(s;t_1,\ldots,t_{r-1},t_r)= \sum_{k=0}^N \sum_{\kappa\vdash k}
c_\kappa(s) p_\kappa(t_1,\ldots,t_{r-1},t_r)
= (s^2-p_1(t_1,\ldots,t_{r-1},t_r))^N,
\label{refnew3}
\eeq
which proves (a).

To prove (b), write $A(s;t_1,\ldots,t_r)$ as 
$A(s;t_1,\ldots,t_r)=sI-B$, where $I$ is the identity
matrix of order $2^r$. Note that $B$ is a matrix with elements
similar to those of $A$ but
with zeros on the diagonal. Introducing a diagonal matrix $D$ of order
$2^r$ by 
\beq
d_{{\pmb l},{\pmb l}}= \left( t_1^{l_1} t_2^{l_2}\cdots t_r^{l_r}
\right)^{1/2},
\eeq
it is easy to see that $DBD^{-1}$ is a real and symmetric matrix.
For such matrices, all eigenvalues are real, with 
geometric multiplicity equal to the algebraic multiplicity. 
Thus this also holds for the eigenvalues of $B$.
The characteristic equation of $B$ is
\beq
\det(B-\la I)=\det(A(\la;t_1,\ldots,t_r))
= (\la^2-\sum_{i=1}^r t_i)^{2^{r-1}}
= (\la-\mu)^{2^{r-1}}(\la+\mu)^{2^{r-1}},
\label{detB}
\eeq
where $\mu=+\sqrt{\sum_{i=1}^r t_i}$.
So for the eigenvalue $\mu$ the geometric
multiplicity is $2^{r-1}$. Since the geometric multiplicity is also
equal to $2^r - \rank(B-\mu I)$, it follows that 
$\rank(B-\mu I)=2^r-2^{r-1}=2^{r-1}$, or $\rank(A(\mu;t_1\ldots,t_r))
=2^{r-1}$, implying the statement (b). 

To prove (c), consider the element $c_{{\pmb l},{\pmb l'}}$ in
the multiplication of $A(s;t_1,\ldots,t_r)$ with 
$A(-s;t_1,\ldots,t_r)$. From the definition (\ref{defA}) it follows
immediately that $c_{{\pmb l},{\pmb l'}}=0$ if ${\pmb l}\ne{\pmb l'}$,
and that $c_{{\pmb l},{\pmb l}}=\sum_i t_i-s^2$.
\mybox

This lemma can now be used to determine the linear (in)dependence
of the weight vectors $X(p;{\pmb k},{\pmb l})$. 
Consider a weight $\la_m=(p-\sum_i m_i,m_1, \ldots, m_n)$, with
all $m_i\geq 0$. Then the multiplicity of $\la_m$ in $\bar V_p$, 
or equivalently the dimension of the weight space $\bar V_p(\la_m)$,
is given by
\beq
d_m=\dim \bar V_p(\la_m) = 2^{\ga(m_1)+\cdots+\ga(m_n)}=2^r,
\eeq
where $r$ is the number of nonzero $m_i$'s.
A basis of $\bar V_p(\la_m)$ is given by the set of vectors
$|p;{\pmb k},{\pmb l}\rangle$ with every $k_i+l_i=m_i$ (or ${\pmb k}+
{\pmb l}={\pmb m}$). 
One can consider another set of $d_m$ vectors $X(p;{\pmb k},{\pmb l})$
with ${\pmb k}+{\pmb l}={\pmb m}$. The coefficient matrix of the
vectors $(-1)^{l_1+\cdots+l_n} X(p;{\pmb k},{\pmb l})$ expressed in
terms of the vectors $|p;{\pmb k},{\pmb l}\rangle$ coincides
with the matrix $A(\sqrt{p};t_1,\ldots,t_r)$ defined in Lemma~3, 
with $t_i$ corresponding to the nonzero $m_i$'s. Thus the determinant
of this matrix is
\beq
(p-\sum_{i=1}^n m_i)^{d_m/2}.
\eeq
In other words, if $\sum_{i=1}^n m_i\ne p$, then the coefficient
matrix is nonsingular, and the $d_m$ vectors $X(p;{\pmb k},{\pmb l})$
with ${\pmb k}+{\pmb l}={\pmb m}$ form a basis for $\bar V_p(\la_m)$.
When $\sum_{i=1}^n m_i=p$, it follows from Lemma~3(b) that the
span of the
$d_m$ vectors $X(p;{\pmb k},{\pmb l})$ with ${\pmb k}+{\pmb l}={\pmb m}$
is a subspace of $\bar V_p(\la_m)$ of dimension $d_m/2$.

\section{The simple module $V_p$}

Denote by $M_p$ the maximal $G=q(n+1)$ submodule of $\bar V_p$ (different
{}from $\bar V_p$ itself). Then the quotient module $V_p=\bar V_p / M_p$
is an irreducible (or simple) $q(n+1)$ module. In this section we shall
show that $V_p$ is finite dimensional, and give its weight
structure, character, and dimension. By definition, $V_p$ and $M_p$
are weight modules, and 
\beq
v\in M_p \Leftrightarrow v_0 \not\in U(G) v.
\label{Mp}
\eeq

For the weight vectors $|p;{\pmb k},{\pmb l}\rangle$ or 
$X(p;{\pmb k},{\pmb l})$ it will be useful to refer to the quantity
$\sum_{i=1}^n (k_i+l_i)=\sum_i m_i$ as the {\em level} of the vector
(or of the corresponding weight).

\begin{prop}
The weight vectors $v$ of $\bar V_p$ satisfy the following~:
\begin{itemize} 
\item[(a)]
if the level of $v$ is greater than $p$ then $v\in M_p$;
\item[(b)]
if the level of $v$ is less than $p$ then $v\not\in M_p$, and denoting
the vectors of the quotient module by their representatives in $\bar V_p$
we can write $v\in V_p$;
\item[(c)]
if the level of $v$ is equal to $p$, consider its weight 
$\la_m=(0,m_1,\ldots,m_n)$. With $d_m=\dim(\bar V_p(\la_m))$, we have
that $\dim M_p(\la_m)=\dim(V_p(\la_m))=d_m/2$. Moreover, 
the vectors $X(p;{\pmb k},{\pmb l})$ of level $p$ with
$k_i+l_i=m_i$ span $M_p(\la_m)$.
\end{itemize}
\end{prop}

\noindent {\bf Proof.}
Consider a fixed weight $\la_m$, and the corresponding weight
vectors $|p;{\pmb k},{\pmb l}\rangle$ and $X(p;{\pmb k},{\pmb l})$.
If the level of $\la_m$ is less than $p$, then the vectors
$X(p;{\pmb k},{\pmb l})$ form a basis for $\bar V_p(\la_m)$ (see end
of section 6). By (\ref{bj-x}) and (\ref{fj-x}),
\[
(b_1^-)^{k_1}(f_1^-)^{l_1}\cdots (b_n^-)^{k_n} (f_n^-)^{l_n} 
X(p;{\pmb k},{\pmb l})
\]
is equal to a nonzero constant times $v_0$. Thus 
$v_0\in U(G) X(p;{\pmb k},{\pmb l})$, in other words all vectors
of weight $\la_m$ are not in $M_p$. 

Let the level of $v$ be greater than $p$. Applying $b_j^-$ or $f_j^-$
reduces the level by one. But for the vectors at level $p+1$ it
follows from (\ref{bj-p}) and (\ref{fj-p}) that the action
of $b_j^-$ or $f_j^-$ yields only vectors of the form 
$X(p;{\pmb k},{\pmb l})$ at level $p$, and all these vectors themselves
are annihilated by $b_j^-$ and $f_j^-$. Thus one deduces
that $v_0$ cannot belong to $U(G)v$ if $v$ has level greater than $p$.

Finally, consider a weight $\la_m=(0,m_1,\ldots,m_n)$ of level $p$.
The vectors $X(p;{\pmb k},{\pmb l})$ with $k_i+l_i=m_i$ are all
annihilated by $b_j^-$ and $f_j^-$, so it follows that they belong
to $M_p$. In this case, we know (see end Section 6) that these
vectors span a subspace of dimension $d_m/2$. Thus $\dim M_p(\la_m)
\geq d_m/2$, and we still need to show that the dimension of
$M_p(\la_m)$ does not exceed $d_m/2$.
The space $\bar V_p(\la_m)$, of dimension $d_m$, is spanned by the
$d_m$ vectors $|p;{\pmb k},{\pmb l}\rangle$ with $k_i+l_i=m_i$.
Assume that $m_n\ne 0$ (the same argument works for another $m_i\ne 0$).
Consider the sets
\bea
S_0 &=& \{ X(p;k_1,l_1,\ldots,k_{n-1},l_{n-1},m_n,0) | k_i+l_i=m_i \},\\
S_1 &=& \{ |p;k_1,l_1,\ldots,k_{n-1},l_{n-1},m_n-1,1\rangle 
 | k_i+l_i=m_i \}, \\
S&=& S_0 \cup S_1.
\eea
Clearly, $\# S_0= \# S_1 = d_m/2$. The vectors in $S_1$ are obviously
linearly independent. By considering the coefficient matrix of the
vectors of $S_0$ in terms of the $d_m$ vectors 
$|p;{\pmb k},{\pmb l}\rangle$, and using Lemma~3, it follows that
the vectors of $S_0$ are also linearly independent, and furthermore
that the vectors of $S$ are linearly independent. Thus $S$ constitutes
a basis for $\bar V_p(\la_m)$. The elements of $S_0$ all belong to 
$M_p(\la_m)$. On the other hand, ${\rm span}(S_1)$ contains no vectors
that are annihilated by $P$. Indeed, consider a linear combination
$c$ of the vectors in $S_1$, and express that $c$ is
annihilated by $b_n^-$. Using (\ref{bj-p}), and linear independence of the
vectors appearing in $b_n^- c$, it follows that $b_n^- c=0$ only if 
all the coefficients in the linear combination $c$ are zero.
One can now deduce that no linear combination of $S_1$ can yield a
vector of $M_p(\la_m)$. This shows that
$\dim M_p(\la_m)=d_m/2$, hence $\dim V_p(\la_m)=d_m/2$. \mybox

It can be verified that $M_p$, which is a $q(n+1)$ module, is 
generated by one vector $w$ as a $q(n+1)$ module (otherwise
said~: $\bar V_p$ containes one $q(n+1)$ highest weight 
singular vector $w$). In our notation, this vector $w$ is equal
to $w=X(p;p,0,\ldots,0, 0,\ldots 0)$.

This proposition gives us explicitly the vectors of $M_p$. Hence
it also gives the (representatives of) the vectors of $V_p=
\bar V_p/M_p$. In particular, $V_p$ is finite dimensional, and
the weight structure of $V_p$ can be deduced. For a weight
$\la_m=(p-\sum_i m_i,m_1,\ldots,m_n)$ (with all $m_i\geq 0$), let
$r$ be the number of nonzero $m_i$'s. Then the multiplicity of
$\la_m$ is $2^r$ if the level of $\la_m$ is less than $p$, and
$2^r/2=2^{r-1}$ if the level of $\la_m$ is equal to $p$ (and,
of course, zero of the level is larger than $p$). 

Once the weight structure is known, it is possible to write down
the character and dimension of $V_p$. To do this, it will be
useful to first determine the decomposition of $V_p$
with respect to the subalgebra $gl(n+1)\subset G$. Since $V_p$
is finite dimensional, it will decompose into a direct sum of
simple finite dimensional $gl(n+1)$ modules. This decomposition
can be derived from the weight structure. The highest weight
is $(p,0,\ldots,0)$, with multiplicity 1. So the $gl(n+1)$ module with
highest weight $(p,0,\ldots,0)$ is a component of the decomposition.
Subtracting the (known) weights of this $gl(n+1)$ module from the 
set of weights of $V_p$, leaves $(p-1,1,0,\ldots,0)$ as the next 
highest weight, also with multiplicity one. Then we go on~:
first subtract all weights (including multiplicities) of the
$gl(n+1)$ module labelled by $(p-1,1,0,\ldots,0)$; then determine the
highest weight of the remaining ones, etc. Finally, one obtains~:
\begin{prop}
The decomposition of the $q(n+1)$ module $V_p$ into $gl(n+1)$ modules 
(with each $gl(n+1)$ module characterized by its highest weight) 
is as follows~:
\beq
V_p \rightarrow (p,0,\ldots,0)\oplus (p-1,1,0,\ldots,0)
\oplus (p-2,1,1,0,\ldots,0)\oplus \cdots \oplus (p-n,1,1,\ldots,1).
\label{decomp}
\eeq
The dimension of $V_p$ is given by~:
\beq
\dim V_p = \sum_{i=0}^n {p-1 \choose i}{p+n-i\choose n-i}.
\label{dimen}
\eeq
\end{prop}
\noindent {\bf Proof.} The decomposition follows from the known
weight structure determined in Proposition~4. From the known
dimension formula (e.g.\ \cite[\S 4]{Wybourne})
of simple $gl(n+1)$ modules, (\ref{dimen}) follows.
\mybox

The (formal) character of a $G$ module $V$ is defined as usual~:
\beq
\ch V = \sum_{\la} \dim V(\la) x_0^{\la_0}\cdots x_n^{\la_n},
\eeq
where $\dim V(\la)$ is the multiplicity of a weight $\la$ in $V$,
and the $x_i$ can be considered as formal variables. 
The character of $V_p$ follows from (\ref{decomp}), using the (known)
characters of $gl(n+1)$ modules. For a $gl(n+1)$ module with
highest weight $\la$, the character is equal to the Schur function
$s_\la(x_0,\ldots,x_n)$. 
In this case, the highest weights appearing in the decomposition
are of a special form; in fact, they are of Frobenius form 
$(p-1-i|i)$~\cite[p.\ 3]{Macdonald}. 
Since in such a case the character
is given by~\cite[p.\ 47]{Macdonald}~:
\beq
s_{(a|b)}= h_{a+1}e_b- h_{a+2}e_{b-1}+ \cdots +(-1)^b h_{a+b+1},
\eeq
where $e_r$ (resp.\ $h_r$) is the $r$th elementary (resp.\ complete)
symmetric function in the $x_i$, it follows that
\beq
\ch V_p = h_{p-n}e_n + h_{p-n-2}e_{n-2}+ \cdots,
\eeq
ending with $h_p$ if $n$ is even and with $h_{p-1}e_1$ if $n$ is odd.

To have the interpretation of $V_p$ as a Fock space of $q(n+1)$,
we still need to show that the Hermitian form is positive definite.

\begin{prop}
The Hermitian form on $V_p$, induced by (\ref{form}), 
is positive definite.
\end{prop}

\noindent {\bf Proof.}
It is clear from (\ref{form}) and (\ref{Mp}) that the Hermitian form is
zero on $M_p$, so (\ref{form}) induces indeed a Hermitian form on $V_p=
\bar V_p/M_p$.
It also follows from (\ref{form}) that $\langle v | w\rangle =0$ if
the weight of $v$ and $w$ is different. So it is sufficient to 
study the behaviour of the Hermitian form on a weight space 
$V_p(\la_m)$ only. 
Let $\la_m$ be fixed, and assume that the level of $\la_m$ is less than
$p$ ($\sum m_i<p$) and that all $m_i$ are nonzero 
(the proof has be slightly changed in the remaining cases
since according to Proposition~4 a different basis must be chosen; 
but in this different basis, it leads essentially to
the same computation).
A basis for $V_p(\la_m)$ is given by the vectors 
$|p;{\pmb k},{\pmb l}\rangle$ with ${\pmb k}+{\pmb l}={\pmb m}$.
Thus we have to show that the matrix $H$ with matrix elements
\beq
H_{{\pmb l},{\pmb l'}}=\langle\ |p;{\pmb k},{\pmb l}\rangle\ |\ 
|p;{\pmb k'},{\pmb l'}\rangle\ \rangle
\eeq
is positive definite. By (\ref{pkkll}), $H_{{\pmb l},{\pmb l'}}$ is equal
to the coefficient of $v_0$ in
\beq
(f_n^-)^{l_n}(b_n^-)^{k_n}\cdots (f_1^-)^{l_1}(b_1^-)^{k_1}
|p;{\pmb k'},{\pmb l'}\rangle.
\eeq
The idea is now as follows~: 
\begin{itemize}
\item The coefficient matrix of the vectors $(-1)^{l_1+\cdots+l_n}
X(p;{\pmb k},{\pmb l})$ expressed in terms of the vectors 
$|p;{\pmb k},{\pmb l}\rangle$ is given by $A=A(\sqrt{p};m_1,\ldots,m_n)$.
Thus the 
coefficient matrix of the vectors $|p;{\pmb k},{\pmb l}\rangle$
expressed in terms of $(-1)^{l_1+\cdots+l_n}
X(p;{\pmb k},{\pmb l})$ is given by $A^{-1}$.
\item The action of 
$(f_n^-)^{l_n}(b_n^-)^{k_n}\cdots (f_1^-)^{l_1}(b_1^-)^{k_1}$
on a vector of the form $(-1)^{l'_1+\cdots+l'_n}
X(p;{\pmb k'},{\pmb l'})$ is diagonal, and determined by
(\ref{bj-x}) and (\ref{fj-x}).
\end{itemize}
This leads to
\beq
H_{{\pmb l},{\pmb l'}}= d({\pmb k},{\pmb l}) 
(A^{-1})_{{\pmb l'},{\pmb l}},
\eeq
where 
\[
d({\pmb k},{\pmb l})=k_1!k_2!\cdots k_n!(p-1)(p-2)\cdots (p-\sum m_i).
\]
It follows that $H=c D^{-1}A^{-T}$, with $c$ a positive constant, 
$A^{-T}$ the transpose of $A^{-1}$, and
$D$ a diagonal matrix with elements 
$D_{{\pmb l},{\pmb l}}=m_1^{l_1}\cdots m_n^{l_n}$.
But $H$ (being symmetric) 
is positive definite if and only if $D^{1/2}HD^{1/2}$ is
positive definite (e.g.\ by the Cholesky decomposition).
Now $D^{1/2}HD^{1/2}=c D^{-1/2} A^{-T} D^{1/2}$; this last matrix
is positive definite if all its eigenvalues are positive.
{}From the proof of Lemma~3(b) (and $A^{-1}$ determined by Lemma~3(c)),
\[
\det(D^{-1/2} A^{-T} D^{1/2}-\lambda I)=
\left((\sqrt{p}-\la)^2-\sum_i m_i \right),
\]
so the eigenvalues are $\la=\sqrt{p}\pm \sqrt{\sum_i m_i}$,
which are indeed positive since $\sum_i m_i<p$. \mybox

\section{Conclusion}

We have given a description of the Lie superalgebra $q(n+1)$
in terms of creation operators $b_i^+$, $f_i^+$ and annihilation
operators $b_i^-$, $f_i^-$ ($i=1,\ldots,n$). The quadratic
relations~(\ref{Q1}) and the triple supercommutation relations
(\ref{Q2}) and (\ref{Q3}) determine the Lie superalgebra 
$sq(n+1)$ completely. The operators $b_i^\pm$ satisfy the relations
of $A$-statistics, and the operators $f_i^\pm$ the relations of
$A$-superstatistics. The combined relations (\ref{Q1})-(\ref{Q3})
can be seen as a unification of $A$-statistics and $A$-superstatistics.

We have shown that $q(n+1)$ has an interesting class
of irreducible representations $V_p$, defined as a quotient module
of an induced module $\bar V_p$. For $p$ a positive integer, these
representations $V_p$ are finite-dimensional, with a unique highest
weight (of multiplicity 1). The Hermitian form that is consistent
with the natural adjoint operation on $q(n+1)$ is shown to be
positive definite on $V_p$. For $q(2)$ these representations
are ``dispin'', since they decompose into the sum of two
irreducible $gl(2)$ representation in the decomposition
$q(2)\supset gl(2)$. Also in the general case, the
decomposition of $V_p$ with respect to $q(n+1)\supset gl(n+1)$
is determined, through the weight structure of $V_p$.
Thus a character and dimension formula for $V_p$ is given.

\section*{Acknowledgements}

T.D.~Palev would like
to thank the University of Ghent for a Visiting Grant, and the
Department of Applied Mathematics and Computer Science for its
kind hospitality during his stay in Ghent. The authors would also
like to thank the referees for their careful reading, for pointing
out some misprints, and for suggesting some improvements.


\begin{thebibliography}{99}

\bibitem{Kac1}
Kac V G 1977
{\em Adv. Math.} {\bf 26} 8 

\bibitem{Kac2}
Kac V G 1978
{\em Lecture Notes in Math.} {\bf 676} 597 

\bibitem{Scheunert}
Scheunert M 1979
{\em The theory of Lie superalgebras}
(Springer, Berlin)

\bibitem{Corwin}
Corwin L, Ne'eman Y and Sternberg S 1975
{\em Rev. Mod. Phys.} {\bf 47} 573 

\bibitem{Balantekin}
Balantekin A B 1984
{\em J. Math. Phys.} {\bf 25} 2028 

\bibitem{Hurni}
Hurni J P 1987
{\em J. Phys. A} {\bf 20} 5755 

\bibitem{VHKT1}
Van der Jeugt J, Hughes J W B, King R C and Thierry-Mieg J 1990
{\em J.\ Math.\ Phys.} {\bf 31} 2278 

\bibitem{VHKT2}
Van der Jeugt J, Hughes J W B, King R C and Thierry-Mieg J 1990
{\em Commun. Algebra} {\bf 18} 3453 

\bibitem{KW}
Kac V G and Wakimoto M 1994
{\em Progress in Math.} {\bf 123} 415 

\bibitem{Serganova1}
Serganova V 1993
{\em Advances in Soviet Math.} {\bf 16} 151 

\bibitem{Serganova2}
Serganova V 1996
{\em Selecta Mathematica} {\bf 2} 607 

\bibitem{VZ}
Van der Jeugt J and Zhang R B 1999
{\em Lett. Math. Phys.} {\bf 47} 49 

\bibitem{pal78}
Palev T D 1987
{\em Funct. Anal. Appl.} {\bf 21} 245 

\bibitem{pal82}
Palev T D 1989
{\em J. Math. Phys.} {\bf 30} 1433 

\bibitem{pal79}
Palev T D 1987
{\em J. Math. Phys.} {\bf 28} 2280 

\bibitem{pal81}
Palev T D 1988
{\em J. Math. Phys.} {\bf 29} 2589

\bibitem{penkov0}
Penkov I 1986
{\em Funct. Anal. Appl.} {\bf 20} 30 

\bibitem{penkov1}
Penkov I and Serganova V 1997 
{\em Lett. Math. Phys.} {\bf 40} 147 

\bibitem{penkov2}
Penkov I and Serganova V 1997
{\em J. Math. Sci.} {\bf 84} 1382 

\bibitem{frappat}
Frappat L and Sciarrino A 1992
{\em J. Math. Phys.} {\bf 33} 3911 

\bibitem{dict}
Frappat L, Sorba P and Sciarrino A 1996
{\em Dictionary on Lie Superalgebras}
(Enslapp-AL-600/96; hep/th/9607161)

\bibitem{palevSA}
Palev T D 1980
{\em J. Math. Phys.} {\bf 21} 1293 

\bibitem{Green}
Green H S 1953 
{\em Phys. Rev.} {\bf 90} 370

\bibitem{Kamefuchi}
Kamefuchi S and Takahashi Y 1960 
{\em Nucl. Phys.} {\bf 36} 177

\bibitem{Ryan}
Ryan C and Sudarshan E C G 1963  
{\em Nucl. Phys.} {\bf 47} 207

\bibitem{Omote}
Omote M, Ohnuki Y and Kamefuchi S 1976 
{\em Prog. Theor. Phys.} {\bf 56} 1948

\bibitem{Ganchev}
Ganchev A Ch and Palev T D 1978 
{\em J. Math. Phys.} {\bf 23} 1100

\bibitem{palevD}
Palev T D 1979 
{\em Czech. Journ. Phys.} {\bf B 29} 91

\bibitem{palevC}
Palev T D 1976 
{\em Lie algebraic aspects of quantum statistics}
(Habilitation Thesis, Sofia)

\bibitem{palevA}
Palev T D 1977
{\em Lie algebraic aspects of the quantum statistics. Unitary quantization
($A$-quantization)}
(Preprint JINR E17-10550; hep-th/9902157)

\bibitem{SNR}
Scheunert M, Nahm W and Rittenberg V 1977
{\em J. Math. Phys.} {\bf 18} 155 

\bibitem{Hughes}
Hughes J W B 1981
{\em J. Math. Phys.} {\bf 22} 245 

\bibitem{Macdonald}
Macdonald I G 1995
{\em Symmetric functions and Hall polynomials}
(Clarendon Press, Oxford)

\bibitem{Wybourne}
Wybourne B G 1970
{\em Symmetry principles and atomic spectroscopy}
(Wiley-Interscience, New York)


\end{thebibliography}
\end{document}